\documentstyle[amscd,righttag,amssymb,intlim]{amsart}
\newcommand{\al}{\alpha}
\newcommand{\bt}{\beta}
\newcommand{\gm}{\gamma}
\newcommand{\Gm}{\Gamma}
\newcommand{\dl}{\delta}
\newcommand{\ve}{\varepsilon}
\newcommand{\sg}{\sigma}
\newcommand{\tht}{\theta}
\newcommand{\om}{\omega}
\newcommand{\Om}{\Omega}
\newcommand{\zt}{\zeta}
\newcommand{\vf}{\varphi}
\newcommand{\nb}{\nabla}
\newcommand{\bC}{\Bbb{C}}
\newcommand{\bN}{\Bbb{N}}
\newcommand{\bR}{\Bbb{R}}
\newcommand{\bZ}{\Bbb{Z}}
\newcommand{\cT}{\cal{T}}
\newcommand{\stto}[1]{\stackrel{#1}{\to}}
\newcommand{\str}{\stackrel}
\newcommand{\pa}{\partial}
\newcommand{\lrc}{\lrcorner}
\newcommand{\sbs}{\subset}
\newcommand{\os}{\overset}
\newcommand{\us}{\underset}
\newcommand{\ol}{\overline}
\newcommand{\wt}{\widetilde}
\newcommand{\wh}{\widehat}
\newcommand{\hra}{\hookrightarrow}
\newcommand{\lmto}{\longmapsto}
\newcommand{\la}{\langle}
\newcommand{\ra}{\rangle}
\newcommand{\arrows}{\aligned & \longrightarrow \\[-3mm] & \longleftarrow
            \endaligned}
\newcommand{\Ext}{\operatorname{Ext}}
\newcommand{\Hom}{\operatorname{Hom}}
\newcommand{\ind}{\operatorname{ind}}
\newcommand{\RRe}{\operatorname{Re}}
\newcommand{\IIm}{\operatorname{Im}}
\newcommand{\sign}{\operatorname{sign}}
\newtheorem{th}{\ \ \ Theorem}[section]
\newtheorem{lem}[th]{\ \ \ Lemma}
\newtheorem{prop}[th]{\ \ \ Proposition}
\newtheorem{cor}[th]{\ \ \ Corollary}

\theoremstyle{definition}
\newtheorem{defi}[th]{\ \ \ Definition}

\theoremstyle{remark}
\newtheorem{rem}[th]{\ \ \ Remark}
\numberwithin{equation}{section}
\begin{document}
%\date{}%
%\dedicatory{}%
%\commby{}%
% ---------------------------------------------------------------
\title[H\"{o}rmander and Maslov Classes, Fomenko's Conjecture]
    {The H\"{o}rmander and Maslov Classes and Fomenko's Conjecture}%
\author{Z. Tevdoradze}%

\address{Faculty of Mechanics and Mathematics I.Javakhishvili Tbilisi state
University\newline
 2,University St., Tbilisi 380043, Republic of Georgia}
\email{tevzaz@@ti.net.ge}%
 \subjclass{58F05, 57R20} \keywords{Symplectic manifold,
connection,
    Maslov class, minimal surface, Hermitian manifold}
%   Abstract.
\begin{abstract}
    Some functorial properties are studied for the H\"{o}rmander classes
    defined for symplectic bundles. The behaviour of the Chern first form
    on a Lagrangian submanifold in an almost Hermitian ma\-ni\-fold is also
    studied, and Fomenko's conjecture about the behaviour of a Maslov
    class on minimal Lagrangian submanifolds is considered.
\end{abstract}
\maketitle
\section*{Introduction}

In this work we are interestied in the characterustic classes of sympletic
bundles and Lagrangian subbundles. One of such classes was discovered
when studying asymptotic solutions of  linear partial differential equations
[1] and called a Maslov class. In [2] Arnold gives a pure geometric
interpretation of the Maslov class.

A generalization of the Maslov class to higher cohomological dimensions
was defined by Arnold and studied in [3]. Another generalization of the
above-mentioned classes is defined by Trofimov in [4].

In [5] and [6] H\"{o}rmander defined the cohomology classes
$\sg(E;X,Y)\in H^1(M;{\bZ})$ for arbitrary sections $X$ and $Y$
of\; Lagrange's\; Grassmanian $L(E)\!\stto{\pi} M$, where $E
\stto{p} M$ is a sympletic fibre bundle. These classes are
discribed in [7].

In [8] Fomenko formulates a conjecture that all Maslov--Arnold
characteristic classes of minimal Lagrangian surfaces are equal to
zero, and this conjecture is proved when Lagrangian surfaces are
submanifolds in $M^{2n}={\bR}^{2n}={\bC}^n$.

The paper is organized as follows. In \S 1 we will describe the condition when
the Maslov index $\ell(\gm)$ is an even number for an arbitrary closed curve
$\gm$ in the Lagrangian manifold (Theorem 1.1) and prove one functorial
property for H\"{o}rmander classes $\sg(E;X,Y)$ (Theorem 1.3).
In \S 2 a class of Lagrangian
manifolds $LH(M)$ is considered in an almost Hermitian manifold $M$ with
the following property: $J^*\al_{{}_H}$ is an exact 1-form in each manifold from
$LH(M)$, where $J$ is an almost complex structure and $H$ is a mean curvature
with respect to the inclusion in $M$. It is proved that $c|_N$ is an exact
2-form for the Lagrangian manifold $N$ in the Hermitian manifold $M$ and
$c|_N=0$, if $N\in LH(M)$, where $c$ is the first Chern form.
The extended variant of Fomenko's conjecture on the class $LH(M)$
also formulated and studied. This conjecture is proved for the class
$LH(T^*Q^n)$, where $Q^n$ satisfies some additional
condition (Theorem 2.8).

In conclusion we note that it remains unknown whether Fomenko's conjecture
is true or not in the general case and how much the class $LH(M)$
extends the class of minimal Lagrangian submanifolds in $M$.

%   1
\section{Construction of Maslov and H\"{o}rmander Classes}

Let $(V^{2n},\om)$ be a symplectic vector space over ${\bR}$. A
real $n$-plane $\zt$ in $V^{2n}$ is called a Lgrangian plane if
the restriction of $\om$ on $\zt$ vanishes. By $L(V)$ we denote
Lagrange's Grassmanian which consists of all Lagrangian
$n$-planes. It is well known that if we have a fixed point $X\in
L(V)$, then $L(V)$ can be represented as a homogeneous space
$U(n)/O(n)$, where the orthogonal group $O(n)$ is cannonically
imbedded in the unitary group $U(n)$.

Let us consider the fibre bundle
$$  SU(n)/SO(n)\to U(n)/O(n)\stto{\det^2} S^1,  $$
where the map $\det^2$ is a fibre map and $SU(n)/SO(n)$ is a fibre
at the point $1\in S^1\sbs {\bC}$. The space $SU(n)/SO(n)$ is
simply connected; therefore the long exact homotopy sequence of
fibre bundles implies that $\pi_1(L(V))={\bZ}$ and hence
$H_1(L(V);{\bZ})={\bZ}$. From the formula of universal
coefficients
$$  0\to \Ext (H_0(L(V),{\bZ}))\to H^1(L(V);{\bZ})\stto{h}
        \Hom(H_1(L(V);{\bZ});{\bZ})\to 0            $$
(where $(h\{f\})\{\sum c_i\otimes g_i\}=\sum f(c_i)\otimes g_i$,
$\{f\}\in H^1(L(V);{\bZ})$, $\{\sum c_i\otimes g_i\}\in
H_1(L(V);{\bZ})$), we obtain $H^1(L(V);{\bZ})=0$.

The differential 1-form $(\det^2)^*\frac{dz}{2\pi iz}$ is the
generator of $H^1(L(V);{\bZ})$, where the form $\frac{dz}{2\pi
iz}$ is the 1-firm on $S^1=\{z\in {\bC}:|z|=1\}\sbs {\bC}$.

%   DEFINITION 1.1
\begin{defi}
A class od cohomology which is defined by 1-form    \linebreak
$(\det^2)^*\frac{dz}{2\pi iz}$ is called a Maslov class.
\end{defi}

%   DEFINITION 1.2
\begin{defi}
A submanifold $N\hra V^{2n}$ is called a Lagrangian manifold if the tangent plane
at each point $x\in N$ is a maximal isotropic plane with respect to
$\om_x=\om|_{{}_{T_xV^{2n}}}$, i.e., $\om_x|_{{}_{T_xN}}=0$.
\end{defi}

For each closed curve $\gm$ in $N$ there is a map
%   (1.1)
\begin{equation}
    G:\gm\to L(V), \quad G(x)=(x,T_xN), \quad x\in \gm,
\end{equation}
which is called Gaussian map. This map defines the integer
%   (1.2)
\begin{equation}
    \ell(\gm)=\int_\gm (\det{}^2\circ G)^*\,\frac{dz}{2\pi iz}
\end{equation}
and thereby the class of cohomology from $H^1(N;{\bZ})$. This
class does not depend on a choice of a Lagrangian plane $X\in
L(V)$ and is called the Maslov class of a Lagrangian manifold $N$.

%   THEOREM 1.3
\begin{th}
If Gaussian map $G:\gm\to L(V)$ can be covered by the continuous map
$\wt{G}:\gm\to U(n)$, i.e., $p\circ \wt{G}=G$, where $p:U(n)\to L(V)$
is the natural projection, then $\ell(\gm)$ is an even number.
\end{th}

%   PROOF
%\begin{pf}
\noindent {\it Proof.\/}
First notice that the map $\det:U(n)\to S^1$ induces an isomorphism
between the fundamental groups $\pi_1(U(n))$ and $\pi_1(S^1)$.
Indeed, the long exact homotopy sequence of fibre bundles
$$  SU(n)\stto{j} U(n) \stto{\det} S^1  $$
contains the sequence
$$  \pi_1(SU(n)) \stto{j_*} \pi_1(U(n)) \stto{(\det)_*}
        \pi_1(S^1) \stto{\pa} \pi_0(SU(n)).     $$
The map $(\det)_*$ is an isomorphism, since $SU(n)$ is a simply
connected space.

It is clear that the map
$$  \pi_1(S^1)\ni [\gm] \to \frac{1}{2\pi i}
        \int_\gm \frac{dz}{z} \in {\bZ}       $$
is an isomorphism. The composition of $(\det)_*$ and the above-mentioned
map give the natural isomorphism
$$  \pi_1(U(n))\ni [\gm] \lmto \frac{1}{2\pi i}
        \int_\gm \frac{d(\det U)}{\det U} \in {\bZ}.  $$
The projection $p:U(n)\to U(n)/O(n)$ defines the monomorpohism
$$  p_* :\pi_1(U(n))\to \pi_1(L(V)),        $$
which is multiplication by 2. This statement immediately follows from
the commutative diagram
$$  \begin{CD}
        U(n) @>p>> U(n)/O(n) \\
        @VV{\det}V @VV{\det^2}V \\
        S^1 @>k>> S^1
    \end{CD},                   $$
where $k(z)=z^2$, $z\in S^1\sbs \bC$. This monomorphism explicitly
can be written as
$$  {\mathbb{Z}}
\ni {\frac{1}{2\pi i}} \int_{\gamma} \frac{d(\det U)}{\det U}
        \str{p_*}{\lmto} \frac{1}{2\pi i}
        \int_{\gamma} \frac{d({\det}^{2} U)}{{\det}^{2} U}\in \bZ,
           $$
where $\gm$ is the closed curve is $U(n)$.

Now from the commutative diagram $p\circ \wt{G}=G$ we have
$\ell(\gm)=[G(\gm)]=p_*[\wt{G}(\gm)]$, where
$[G(\gm)]\in \pi_1(L(V))$ and $[\wt{G}(\gm)]\in \pi_1(U(n))$. Clearly,
$$  \frac{\ell(\gm)}{2}=\frac{1}{2\pi i}\int_{\wt{G}(\gm)}
        \frac{d(\det U)}{\det U}\in \bZ. \;\;\qed       $$
%\end{pf}

%   REMARK 1.4
\begin{rem}
The experience accumulated in investigating various concrete mechanical
systems and calculations of the Maslov class for
the Liouville tori of such systems show that an overwhelming
majority of systems have got even numbers as components of
the Maslov class [9], [10]. Theorem 1.3 describes
the mechanism of such occurrence.
\end{rem}

Now we briefly  recall the definition and construction of more
general classes, namely, of H\"{o}rmander classes [7].

Let $X$, $Y$, $Z$, $W$ be four Lagrangian $n$-planes in $V^{2n}$ such that
$Z\cap X=Z\cap Y=\{0\}$ and $W\cap X=W\cap Y=\{0\}$. For them
H\"{o}rmander  defined the index at the intersection
%   (1.3)
\begin{equation}
    (X,Y,Z,W)=\frac{1}{2}(\sign Q_Z-\sign Q_Z)=
        \ind Q_W-\ind Q_W,
\end{equation}
where $Q_Z$ and $Q_W$ are nondegenerate quadratic forms on
$Y/X\cap Y$ defined by the formulas
$$  Q_Z(y)=\om(p_Z^Xy,y), \quad Q_W(y)=\om(p_W^X(y),y)  $$
(here $p_Z^X$ and $p_W^X$ are the projections of $Y/Y\cap X$ onto $Z$ and $W$,
respectively, accross $X$). The following relations immediately follow
from (1.3):
%   (1.4)
\begin{equation}
\begin{aligned}
    & (X,Y,Z,W)=-(X,Y,W,Z),  \\
    & (X,Y,Z,W)+(X,Y,W,V)+(X,Y,V,Z)=0,  \\
    & (X,Y,Z,W)=-(Z,W,X,Y).
\end{aligned}
\end{equation}

Let $E \stto{p} M$ be a symplectic vector fibre bundle, i.e., for
each point $m\in M$ there is a symplectic 2-form $\om_m$ on $E_m$
which smoothly depends on $m$. Then we can consider the fibre
bundle $L(E) \stto{\pi} M$, where $L(E)_m=L(E_m)$. Now let $X$ and
$Y$ be the sections of the bundle $L(E)$, and $U=\{U_\al\}_{\al\in
I}$ be an open covering of $M$ such that all nonempty finite
intersections are diffeomorphic to ${\bR}^n$ (such a covering is
called a good covering). Then for arbitrary neighbouhoods $U_\al$
and $U_\bt$ one can find sections $Z_\al$ and $Z_\bt$ of $L(E)$ on
$U_\al$ and $U_\bt$, respectively, such that $Z_\al(m)$ and
$Z_\bt(m)$ are transversal to $X(m)$ and $Y(m)$ for each point
$m\in U_\al\cap U_\bt$. In this situation the \v{C}ech 1-cocycle
$\sg$ can be defined by the formula
%   (1.5)
\begin{equation}
    \sg(U_\al,U_\bt)=(X,Y,Z_\al,Z_\bt).
\end{equation}

The cohomology class $[\sg]$ is called the H\"{o}rmander class for
a symplectic vector fibre bundle $E$ and sections $X$, $Y$ and
denoted by $\sg(E;X,Y)\in \check{H}^1(M;{\bZ})$.

From (1.4) and (1.5) one can conclude that
%   (1.6)
\begin{equation}
    \sg(E;X,Y)=-\sg(E;Y,X).
\end{equation}

Let $(E',p',L(E))$ be a pullback bundle of the bundle $E$ to $L(E)$:
$$  \begin{array}{ccc}
        E' & \longrightarrow & E \\[3mm]
        \Big\downarrow\vcenter{\rlap{$\scriptstyle{p'}$}}
            &&\Big\downarrow\vcenter{\rlap{$\scriptstyle{p}$}} \\[3mm]
        L(E) & \os{\pi}{\longrightarrow} & M
    \end{array}
        \ \hskip+1cm
    \begin{array}{ccc}
        L(E') & \us{s'}{\os{\pi'}{\arrows}} & L(E) \\
        \vcenter{\llap{$\scriptstyle{\cT}$}} \Big\uparrow \!
            \Big\downarrow\vcenter{\rlap{$\scriptstyle{\pi'}$}}
            &&\Big\downarrow\vcenter{\rlap{$\scriptstyle{\pi}$}}\\
        L(E) & \us{s}{\os{\pi}{\arrows}} & M
    \end{array}.            $$
Then we can consider the fibre bundle $L(E) \stto{\pi'} L(E)$ and
the natural section $\cT:L(E)\to L(E')$ which is defined by the formula
$\cT(x)=x$, where $x\in L(E)_m$ and the fibres $E_x'$, $E_m$ are identified.
An arbitrary section $S:M\to L(E)$ can be lifted to the section
$S':L(E)\to L(E')$ by the formula $S'(x)=S(m)$, where $x\in L(E)_m$,
$m\in M$, i.e., if $x=(m,\xi)$ and $S(m)=(m,\eta)$, then $S'(x)=(m,\xi,\eta)$
(here, as before, $L(E')_x$ and $L(E)_m$ are identified).

One can verify that $Y^*(S')=S$ and $Y^*(\cT)=Y$ for each section $Y$
of the bundle $L(E)$.

It is not difficult to see that class $[\sg]$ has got the following
functorial property: if $(f^*E,f^*(p),M_1)$ is a pullback bundle
of the bundle $(E,p,M)$, where $f:M_1\to M$ is a smooth map, then
%   (1.7)
\begin{equation}
    f^*\sg(E;X,Y)=\sg(f^*E; f^*X,f^*Y)
\end{equation}
$(f^*X$, $f^*Y$ are pullback sections of the sections $X$ and $Y$
by the map $f$).

The following equalities arise from (1.4) and (1.6):
%   (1.8)
\begin{equation}
\begin{aligned}
    & \sg(L(E');X',\cT)-\sg(L(E');Y',\cT)=  \\
    & \ \hskip+2cm =\sg(L(E');X',Y')\in
        \check{H}^1(L(E);\bZ),  \\
    & \sg(E;X,Y)=Y^*\sg(L(E');X',\cT).
\end{aligned}
\end{equation}

For each section $X$ of the bundle $L(E)$ there is a characteristic class
$\sg_X$ defined by the formula
%   (1.9)
\begin{equation}
    {\sg_X \equiv \sg(L(E);X',\cT)}\in
    \check{H}^{1}
    (L(E);{\bZ}),
\end{equation}
such that

\begin{enumerate}
\item[i)] the restriction of $\sg_X$ on a fibre $L(E)_m$, $m\in M$,
    is the generator of $\check{H}^1(L(E)_m;\bZ)$;
\item[ii)] $X^*\sg_X=0$.
\end{enumerate}

The property ii) immediately follows from (1.8) and (1.6). Indeed, \linebreak
$X^*\sg_X=X^*\sg(L(E');X',\cT)=\sg(E;X^*(X'),X^*(\cT))=\sg(E;X,X)=0$.
The first property is more difficult and can be found in [7].

To prove Theorem 1.6 we will use the following well-known theorem [11].

%   THEOREM 1.5
\begin{th}[Dold--Thom--Gysin]
Let $h^*$ be a multiplicative cohomology theory and let
$F\stto{i} E \stto{\pi} B$ be an $h^*$-fibration. Suppose there are
elements $a_1,\dots,a_r$ in $h^*(E)$ such that $(i^*a_1,\dots,i^*a_r)$
is an $h^*($point$)$-base for $h^*(F)$ as an $h^*($point$)$-module; then
$(a_1,\dots,a_r)$ is an $h^*(B)$-base for $h^*(B)$-base for $h^*(E)$
as an $h^*(B)$-module.
\end{th}

%   THEOREM 1.6
\begin{th}
If $X$ and $Y$ are sections of the bundle $L(E)$, and $X'$, $Y'$
are the lifted sections of the bundle $L(E')$, then
%   (1.10)
\begin{equation}
    \sg(L(E);X',Y')=\pi^*\sg(E;X,Y).
\end{equation}
\end{th}

%   PROOF
\begin{pf}
By virtue of the properties (i) and (ii) of the H\"{o}rmander class $\sg_X$ and
Theorem 1.5 we can write
%   (1.11)
\begin{equation}
    \sg_Y=\pi^*[\bt(X,Y)]+k(X,Y)\sg_X,
\end{equation}
where $\bt(X,Y)$ is the closed 1-form on $M$ and $k(X,Y)\in
{\bZ}$. Applying $X^*$ and $Y^*$ to both sides of (1.11) and
taking into account (1.8) and (ii), we respectively have
%   (1.12)-(1.13)
\begin{align}
    & \sg(E;X,Y)=[\bt(X,Y)]+k(X,Y)\cdot 0=[\bt(X,Y)],  \\
    & Y^*\sg_Y=[\bt(X,Y)]+k(X,Y)\sg(E;X,Y)=0.
\end{align}
Now from (1.6), (1.12), (1.13) we have obtain $k(X,Y)=1$, and formula (1.11)
acquires the form
%   (1.14)
\begin{equation}
    \sg_Y=\pi^*(E;Y,X)+\sg_X.
\end{equation}
It is clear that (1.14) and (1.10) are equivalent equalities (see (1.8)).
\end{pf}

%   2
\section{Fomenko's Conjecture}

Let $M^{2n}$ be an almost Hermitian manifold, i.e., there exist
a tensor field $J$ and a Riemannian metric $g$ on $M$, where $J$
is an endomorphism of the tangent bundle $TM$ with the property
$J^2=-id$ and $g$ is an invariant with respect to $J$, i.e.,
$g(JX,JY)=g(X,Y)$ for all vector fields $X,Y\in \cT(M)$ on $M$.
The fundamental 2-form on $M$ is defined by the formula
%   (2.1)
\begin{equation}
    \Phi(X,Y)=g(X,JY), \quad X,Y\in \cT(M).
\end{equation}
When $\Phi$ is a closed 2-form, $M$ is called an almost Kaehlerian manifold
and $\Phi$ is called an almost Kaehler form on $M$. In this case $M$
is a symplectic manifold with a symplectic 2-form $\Phi$. The converse
is also true: for every symplectic manifold $(M^{2n},\om)$ there exists
an almost Hermitian structure on $M$ such that the fundamental 2-form
of this structure coincides with the symplectic 2-form $\om$.
Thus every symplectic manifold is an almost Kaehlerian manifold [12].

Now let $N^n$ be a Lagrangian submanifold in an almost Kaehlerian manifold
$M^{2n}$, i.e., $i^*\Phi=0$, where $i:N\to M$ is the imbedding map.
By $VN$ we denote the normal fibre bundle on $N$ with respect to the metric $g$
on $M$ (i.e., $V_xN$ is a $n$-plane in $T_xM$, which is orthogonal
to $T_xN$, $x\in N$), and by $\nb$ we denote the Riemannian connection
on $M$. We recall that the trace of the following morphism of the fibrations
$$  B:TN \oplus TN \to VN, \quad B(x,y)=(\nb_xY)^V      $$
(where $x,y\in T_nN$, $n\in N$, and $Y$ is the locally defined vector field on
$M$ which extends the vector $y$, $Y|_N$ is the section of $TN$, and
$(\nb_xY)^V$ is the normal component of vector $\nb_xY$), which can also be
considered as a section of the fibration $VN$, is called the mean
curvative of the submanifold $N$ and denoted by $H$. More precisely,
if $e_1,\dots,e_n$ is an orthobasis of the plane $T_nN$ and $V$ is a locally
defined vector field on $M$ whose  restriction on $N$ is a section of $VN$ and
$V(n)=v$, then the mean curvature $H$ is defined by the relation
%   (2.2)
\begin{equation}
    g(H(n),v)=-\sum_{i=1}^n g(\nb_{e_i}V,e_i).
\end{equation}

We recall that $N\str{i}{\hra} M$ is called a locally minimal submanifold
if the mean curvature $H$ is zero at every point of $N$. From (2.1)
we can conclude that the diagram
\vskip+0.5cm

\noindent is commutative. By $\al_{{}_H}$ we denote
a section of the fibration $A^1(VN)$
(the fibration of the exterior 1-forms on $VN$) defined as
%   (2.3)
\begin{equation}
    \al_{{}_H}(n)(X)=g(H(n),X), \quad n\in N, \quad X\in V_nN.
\end{equation}

%   DEFINITION 2.1
\begin{defi}
We say that the Lagrangian manifold $N$ in $M$ is locally minimal
from the cohomological point of view if the 1-form $\bt=J^*\al_{{}_H}$
defines the trivial class of cohomology $H^1(N;\bR)$.
\end{defi}

The class of the above-defined submanifolds is denoted by $LH(M)$.

In [13] Le Hong Van, relying on the calibrated geometry
methods developed in the fundamental work [14], shows that there exists
an 1-form $\psi$ on $L(TM)$ such that the Lagrangian submanifold
$N$ in $M$ is locally minimal if and only if
$$  G^*(\psi)=0,        $$
where $G:N\to L(TM)$ is the Gaussian map defined as $G(x)=(x,T_xN)\in L(TM)_x$,
$x\in N$.

Now we will briefly describe the 1-form $\psi$.

Let $T_x^c(M^{2n})$ be the complexification of the tangent space of $M$
at point $x$, $x\in M^{2n}$. By $S_x^c$ we denote a complex  subspace of
$T_x^c(M^{2n})$ which contains all eigenvectors of the operator $J$
with the eigenvalue $i$. The unitary bases $R^c(x)=\{\ve_1(x),\dots,\ve_n(x)\}$
of the space $S_n^c(x)$ sonstruct the principal bundle $U(M)$ over $M$,
with the strucrure group $U(n)$. By $\ve_{\ol{1}},\dots,\ve_{\ol{n}}$
we denote the complex conjugate vectors of the vectors
$\ve_1(x),\dots,\ve_n(x)$. Then the vectors
\begin{align*}
    e_i(x) & =\frac{1}{\sqrt{2}}(\ve_i(x)+\ve_{\ol{i}}(x)),
        \quad i=\ol{1,n},  \\
    e_{\ol{i}}(x) & =\frac{i}{\sqrt{2}}(\ve_i(x)-\ve_{\ol{i}}(x))=Je_i,
        \quad i=\ol{1,n},
\end{align*}
make up families of orthonormal real vectors.

An infinitesimal connection on the principal bundle $U(M)$ is called
an almost Hermitian connection on $M$.

Let $U=\{U_\al\}_{\al\in I}$ be the open covering of the manifold $M$, which
is equipped with local sections $S_\al:U_\al\to U(M)$; then an almost
Hermitian connection $\pi$ can be defined by the given 1-forms $\pi_\al$
on $U_\al$, $\al\in I$, with values in the Lie algebra $u(n)$.
We can express $\pi_\al$ in terms of the $(\pi_j^i)_{i,j=\ol{1,n}}$
matrix, where $\pi_j^i$ are the 1-forms on $U_\al$ and
$\pi_j^i+\ol{\pi}_i^j=0$. If we define $S_\al$ as $S_\al(x)=R_\al^c(x)$,
$x\in U_\al$, then there exist metrices $B_\bt^\al(x)\in U(n)$ such that
\begin{gather*}
    S_\al(x)B_\bt^\al(x)=S_\bt(x),  \\
    \pi_\bt=(B_\bt^\al)^{-1}\pi_\al B_\bt^\al+
        (B_\bt^\al)^{-1}dB_\bt^\al, \quad \al,\bt\in I, \;\;
        x\in U_\al\cap U_\bt.
\end{gather*}

Frames of the type $R^c(x)$ and
$\ol{R}{}^c(x)=\{\ve_{\ol{1}}(x),\dots,\ve_{\ol{n}}(x)\}$ construct
the principal bundle $U'(M)$ on $M$ with a structure group $U(n)$ which is
a subgroup in $U(2n)$ by means of the imbedding
$A\to \begin{pmatrix} A & 0 \\ 0 & \ol{A} \end{pmatrix}$, $A\in U(n)$.
Then the above-mentioned Riemannian connection $\nb$ with respect
to the metric $g$ can be expressed on each open set $U_\al$ by the matrix
of 1-forms
$$  \begin{pmatrix}
        \pi_j^i & \pi_{\ol{j}}^i \\[2mm]
        \pi_j^{\ol{i}} & \pi_{\ol{j}}^{\ol{i}}
    \end{pmatrix},      $$
where $\pi_{\ol{j}}^{\ol{i}}=\ol{\pi_j^i}$,
$\pi_j^{\ol{i}}=\ol{\pi_{\ol{j}}^i}$ (since $\nb$ is a real connection) and
$\pi_j^{\ol{i}}+\pi_{\ol{j}}^i=0$ (for $\nb g=0$).

The matrices $(\pi_j^i)$ (on each $U_\al$) define an almost Hermitian
connection which is called the first canonical connection for an almost
Hermitian ma\-ni\-fold $M$.

From the first structure equations of E. Cartan we have
%   (2.4)
\begin{equation}
    d\tht^i=\tht^k\wedge \pi_k^i+\tht^{\ol{k}}\wedge \pi_{\ol{k}}^i,
\end{equation}
where $\tht^j=\frac{1}{\sqrt{2}}(e_j)^*+i(Je_j)^*$, $j=\ol{1,n}$,
$k=\ol{1,n}$. (2.4) is equivalent to
$$  d\tht^i=\tht^k\wedge \pi_k^i+\Gm_{\ol{k}q}^i\tht^{\ol{k}}
        \wedge \tht^q+\frac{1}{2}\,T_{\ol{k}\ol{q}}^i\tht^{\ol{k}}
        \wedge \tht^{\ol{q}}, \quad i=\ol{1,n}, \;\;
        k=\ol{1,n},\;\; q=\ol{1,n},     $$
where $\Gm_{\ol{k}q}^i$ are the coefficients of the Riemannian connection with
respect to frames from the space $U'(M)$, and
$T_{\ol{k}\ol{q}}^i=\Gm_{\ol{k}\ol{q}}^i-\Gm_{\ol{q}\ol{k}}^i$.
Now we note that the forms
$T_{\ol{k}\ol{q}}^i \tht^{\ol{k}}\wedge \tht^{\ol{q}}$, $i=\ol{1,n}$,
construct the vector 2-form on $U(M)$ which coincides with the torsion
vector $T$ of an almost cemplex structure. By direct calculation
(using the fact that $d\Phi=0$) we can obtain $\Gm_{\ol{k}q}^i=0$. So (2.4)
finally has the following view
%   (2.5)
\begin{equation}
    d\tht^i=\tht^k\wedge \pi_k^i+\frac{1}{2}\,T_{\ol{k}\ol{q}}^i
        \tht^{\ol{k}}
        \wedge \tht^{\ol{q}}, \quad i=\ol{1,n}, \;\;
        k=\ol{1,n},\;\; q=\ol{1,n}.
\end{equation}

By $L(TM)$ we denote Lagrange's Grassmanian associated with
the symplectic fibration $TM\stto{p} M$. The map $q:U(M)\to L(TM)$,
defined as $q((\ve_1,\dots,\ve_n))=e_1\wedge \cdots \wedge e_n$
(here the multivector $e_1\wedge \cdots \wedge e_n$ is identified
with the Lagrange $n$-plane span $\{e_1,\dots,e_n\}$) defines
the principal bundle with the structure group $O(n)$. On the space
$U(M)$ there exists an 1-form
%   (2.6)
\begin{equation}
    \ol{\psi}=-\bigg(\sum_{k=1}^n i\pi_k^k +2\IIm\bigg(
        \sum_{i,k=1}^n T_{\ol{i}\ol{k}}^i \tht^{\ol{k}}\bigg)\bigg),
\end{equation}
which can be expressed as a pullback form by the map $q$. Indeed,
the 1-form $\sum\limits_k\pi_k^k$ vanishes on the fibres of
$q:U(M)\to L(TM)$ and is invariant under the action of the group $O(n)$ on
$U(M)$. Therefore the form $\sum\limits_ki\pi_k^k$ can be pulled down
on the $L(TM)$. It is easy to verify that
$$  \sum_i(\ve_{\ol{i}}\lrc T^i)=\sum_{i,k}T_{\ol{i}\ol{k}}^i
        \tht^{\ol{k}} -\sum_{i,k} T_{\ol{k}\,\ol{i}}^i \tht^{\ol{k}}=
        2\sum_{i,k} T_{\ol{i}\,\ol{k}}^i \tht^{\ol{k}},     $$
and if $R_\al(x)=R_\bt(x)A_\al^\bt(x)$, where
$A_\al^\bt(x)=(A_j^{i'})_{i',j=\ol{1,n}}\in O(n)$, $\al,\bt\in I$,
then
$$  \sum_i(\ve_{\ol{i}}\lrc T^i)=\sum_{i,i'}A_i^{i'}
        (\ve_{\ol{i}'} \lrc T^i)=\sum_{i,i',p} A_i^{i'}
        (\ve_{\ol{i}'}\lrc T^{p'})A_{p'}^i=
        \sum_{i'}(\ve_{\ol{i'}}\lrc T^{i'}).        $$
This means that the form $2\IIm(T_{\ol{i}\,\ol{k}}^i \tht^{\ol{k}})$
can also be pulled down and so there exists a 1-form $\psi$ on $L(TM)$
such that
%   (2.7)
\begin{equation}
    \ol{\psi}=q^*\psi.
\end{equation}

When $M$ is a Hermitian manifold $(T=0)$, the 1-form $\psi$ on
$L(TM)$ has a simpler form [8]
%   (2.8)
\begin{equation}
    \psi=J\,df+d\tht,
\end{equation}
where $f$ and $\tht$ are the locally defined functions on $L(TM)$.

%   LEMMA 2.2
\begin{lem}
The diagram is commutative:
\vskip+0.5cm

\noindent where $i$ is the imbedding map, $G$ is the Gaussian map.
\end{lem}

%   PROOF
%\begin{pf}
\noindent {\it Proof.\/}
Let $x\in N$, $D$ be an open neighbourhood of $x$ in $N$ and $S_D$ be
a local section $S_D:D\to U(M)$, $S_D(x)=(\ve_1(x),\dots,\ve_n(x))$.
Let $S_{D'}$ be the extension of $S_D$ to the tubular neighbourhood
$D'$ of $D$ in $M$; then a locally $n$-form on $D'$
$$  \vf=\RRe(S_{D'}^* (\tht^1\wedge \tht^2\wedge \cdots
        \wedge \tht^n)),        $$
is a $n$-form of comass 1 on $D'$ and $\vf|_D\equiv 1$. If $X\in V_nN$,
$n\in N$, we have
\begin{align*}
    (X\lrc d\vf)(e_1,\dots,e_n) & =\sum_{i=1}^n (-1)^i
        e_i(\vf(X,e_1,\dots,\wh{e}_i,\dots,e_n))+  \\
    & +X(\vf(e_1,\dots,e_n))+  \\
    & +\sum_{i<j} (-1)^{i+j} \vf([e_i,e_j],X,\dots,\wh{e}_i,\dots,
        \wh{e}_j,\dots,e_n)+  \\
    & +\sum_{i=1}^n (-1)^i \vf([X,e_i],e_1,\dots,\wh{e}_i,\dots,e_n).
\end{align*}
Since $\vf$ has comass 1, $X(\vf(e_1,\dots,e_n))=0$, and $X\lrc \vf\equiv 0$
for all $X\bot T_nN$. Thus from the above equality we have
\begin{align*}
    (X\lrc d\vf)(e_1,\dots,e_n) & =\sum_{i=1}^n (-1)^i
        \vf([X,e_i],e_1,\dots,\wh{e}_i,\dots,e_n)=  \\
    & =-\sum_i g([X,e_i],e_i)=-\sum_i g(\nb_Xe_i'-\nb_{e_i}X,e_i)=  \\
    & =\sum_i g(\nb_{e_i}X,e_i)=-g(H,X).
\end{align*}
In the fourth equality the fact is used that $\nb_X\la e_i',e_i'\ra=0$,
and $e_i'$ denotes the locally defined vector fields on $D'$ which extend
the vectors $e_i$, $i=\ol{1,n}$. By direct calculation we can obtain
$(S_D)^*(\ol{\psi})(JX)=(X\lrc d\vf)(e_1,\dots,e_n)=-\al_{{}_H}(X)$.
Since $q\circ S_D=G$, we have
$$  \al_{{}_H}(X)=-G^*\psi(JX)=-J^*G^*\psi(X).      $$
Therefore
%   (2.9)
\begin{equation}
    \bt=J^*\al_{{}_H}=G^*\psi.  \;\;\qed
\end{equation}
%\end{pf}

The second structure equations of E. Cartan for the first canonical
connection $\pi=(\pi_j^i)_{i,j=\ol{1,n}}$ have the form
%   (2.10)
\begin{equation}
    \Om=d\pi+\pi\wedge \pi,
\end{equation}
where $\Om$ is the curvature form for the connection $\pi$. Taking into
account $\pi_j^i=-\ol{\pi}_i^j$, from (2.10) we have
$\sum\limits_k\Om_k^k=-\sum\limits_k\ol{\Om}_k^k$ and
$d(\sum\limits_k \pi_k^k)=\sum\limits_k\Om_k^k$. The form
$\frac{i}{2\pi} \sum\limits_k \Om_k^k$ is real and can be pulled down on $M$.
It is called the Chern first form and denoted by $c$. By (2.6) we obtain
%   (2.11)
\begin{equation}
    d\ol{\psi}=-2\pi p^*c-2d\bigg(\IIm\sum_{i,k}
        T_{\ol{i}\,\ol{k}}^i\tht^{\ol{k}}\bigg).
\end{equation}

For an open covering $U=\{U_\al\}_{\al\in I}$ of the manifold $M$ there
exists a commutative diagram
\vskip+0.5cm

\noindent Since $p^*c=q^*\pi^*c$ and $G^*\pi^*=i^*$, we have
$$  S_{U_\al}^*d\ol{\psi}=-2\pi c|_{N_\al}-2S_{U_\al}^*
        d\bigg(\IIm \sum T_{\ol{i}\,\ol{k}}^i
        \tht^{\ol{k}}\bigg).                $$
Now, taking into account equality (2.9) and $S_{U_\al}^* q^*=i_\al^*G^*$,
we obtain
%   (2.12)
\begin{equation}
    i_\al^*d\bt=i_\al^* G^*d\psi=-2\pi c|_{N_\al}-2d\bigg(S_{U_\al}^*
        \IIm \bigg(\sum_{i,k} T_{\ol{i}\,\ol{k}}^i
        \tht^{\ol{k}}\bigg)\bigg).
\end{equation}
Thus $\bt$ is the global 1-form on $N$, the expression
$d(S_{U_\al}^*\IIm (\sum\limits_{i,k} T_{\ol{i}\,\ol{k}}^i \tht^{\ol{k}})$
defines the global closed 2-form on the Lagrangian submanifold
$N$ and we denote it by $\gm$. We have

%   PROPOSITION 2.3
\begin{prop}
The Chern form $c$ on the Lagrangian submanifold $N$ in an almost
Hermitian manifold $M$ can be expressed by the sum
%   (2.13)
\begin{equation}
    i^*c=-\frac{1}{2\pi}\,d\bt-\frac{1}{\pi}\,\gm,
\end{equation}
where $i$ is the imbedding map $i:N\to M$, the forms $\bt$ and $\gm$
depend on the mean curvature and torsion vector form, respectively.
\end{prop}

%   THEOREM 2.4
\begin{th}
If $M$ is a Hermitian manifold then:

\rom{i)} the form $\pi^* c$ is an exact form on $L(TM)$;

\rom{ii)} if $N$ is a Lagrangian submanifold in $M$ then
    $i^*c=-\frac{1}{2\pi}\,d\bt$;

\rom{iii)} the form $\psi$ is a closed form exactly when $c=0$;

\rom{vi)} if $\bt$ is closed on $N$ $($in particular, when $N\in LH(M))$
    then $i^*c=0$;

\rom{v)} if $M$ is the Kaehlerian manifold then $\psi$ is closed exactly
    when Ricchi's tensor is identicaly equal to zero.
\end{th}

%   PROOF
\begin{pf}
From formulas (2.11), (2.9) we have
$$  q^* d\psi =-2\pi p^* c.     $$
Taking into account $\pi\circ q=p$, we obtain
$$  q^*(d\psi+2\pi\cdot \pi^* c)=0.     $$
As $q:U(M)\to L(TM)$ is a fibration, the above equation implies
%   (2.14)
\begin{equation}
    \pi^* c=-\frac{1}{2\pi}\,d\psi.
\end{equation}
By formula (2.14) it is clear that $\psi$ is closed exactly when $c=0$.
Now we recall that the form $c$ on the Kaehlerian manifold can be
locally expressed by the formula (see [12])
%   (2.15)
\begin{equation}
    c=\frac{1}{2\pi} \sum_{i,j} R_{i\ol{j}}dz^i\wedge dz^{\ol{j}}.
\end{equation}
Formulae (2.14), (2.15) prove assertions i), iii), v). (2.13) immediately
implies ii) and iv).
\end{pf}

When $M={\bC P}^{n}$ and $N\hra {\bC P}^{n}$ is Lagrangian
manifold, the conditions of Theorem 2.4 are fulfilled for Chern
forms of higher dimension on ${\bC P}^{n}$. This result follows
from the well-known fact that the Chern forms on ${\bC P}^n$ are
calculated by the formula
$$  c_k=m(k) c^k, \quad k=\ol{1,n}, \quad m(k)\in\bN.       $$
Now lat $X$ and $Y$ be the sections of the bundle $i^*L(TM)$ on $N$,
where $N$ is a submanifold in $M$. By Theorem 1.5, for every closed 1-form
$\vf$ on $L(TM)$ there exists a 1-form $\dl$ on $N$ such that
%   (2.16)
\begin{equation}
    \wh{i}^*[\vf]=\pi_1^*[\dl]+k([\vf],X)\sg_X, \quad
        k([\vf],X)\in \bR;
\end{equation}
here $\wh{i}$ is an imbedding map $i^*L(TM)\stto{\wh{i}} L(TM)$ and
$\pi_1$ is the projection of the bundle $i^*L(TM)$. Applying $X^*$
and $Y^*$ to (2.16), we obtain respectively
%   (2.17)
\begin{equation}
    X^*\wh{i}^* [\vf]=[\dl]; \quad Y^*\wh{i}^*[\vf]=[\dl]+
        k([\vf],X)\sg(i^*(TM);X,Y).
\end{equation}
From the above equations we have
%   (2.18)
\begin{equation}
    k([\vf],X)\sg(i^*(TM);X,Y)=Y^*\wh{i}^*[\vf]-X^*\wh{i}^*[\vf].
\end{equation}
In the particular case, when $\wh{i}^*[\vf]$ is not a pullback class
from $H^1(N;\bR)$, we have $k\neq 0$ and therefore (2.18) can be rewritten as
%   (2.19)
\begin{equation}
    \sg(i^*(TM);X,Y)=\frac{1}{k([\vf],X)}\,(Y^*\wh{i}^*[\vf]-
        X^*\wh{i}^*[\vf]).
\end{equation}
Thus we have proved

%   PROPOSITION 2.5
\begin{prop}
For an arbitrary closed $1$-form $\vf$ on $L(TM)$ (it is assumed that
$\wh{i}^*[\vf]$ is not a pullback class from $H^1(N;\bR)$ by the map $\pi_1)$,
the H\"{o}rmander class $\sg(i^*(TM);X,Y)$ can be expressed
by formula $(2.19)$.
\end{prop}

%   PROPOSITION 2.6
\begin{prop}
If $N$ is a manifold from $LH(M)$ and $\psi$ is the closed $1$-form
on $L(TM)$, then
%   (2.20)
\begin{equation}
    \wh{i}^*[\psi]=k([\psi],G)\sg_G.
\end{equation}
\end{prop}

%   PROOF
\begin{pf}
For the 1-form $\psi$ and the Gaussian section $G$ (2.16) can be
rewritten as
$$  \wh{i}^*[\psi]=\pi_1^*[\dl]+k([\psi];G)\sg_G.       $$
Applying $G^*$ to the formula and taking into account $N\in LH(M)$ and
Lemma 2.2, we conclude that $[\dl]=0$.
\end{pf}

Now let $M^{2n}=T^*Q^n$ be a cotangent bundle and $\om$ be a canonical
symplectic structure on $M^{2n}$. $M^{2n}$ can be considered as an almost
Kaehlerian manifold with the fundamental 2-form $\om$. There is a fixed
section $A$ of the fibre bundle $L(TM^{2n}) \stto{\pi} M^{2n}$,
which is defined at
every point $m\in M^{2n}$ as a tangent plane at the point $m$ of the fibre
$p^{-1}(p(m))$, where $p$ is a natural projection $T^*Q^n\stto{p} Q^n$.

For every Lagrangian submanifold $N^n\str{i}{\hra} M^{2n}$, the Gaussian map
$G:N^n\to i^* L(TM^{2n})$ is defined as previously. The following definition
is equivalent to Definition 1.1.

%   DEFINITION 2.7
\begin{defi}
The cohomology class $\sg(i^*(TM);G,A)\in H^1(N;\bZ)$ is called
the Maslov class of the Lagrangian submanifold $N$ and we denote
it by $\sg_N$.
\end{defi}

%   THEOREM 2.8
\begin{th}
Let $N^n$ be a manifold from $LH(M^{2n})$. Then:
\begin{enumerate}
\item[i)] if $H^1(Q^n;{\bR})=0$ and $\psi$ is a closed $1$-form, then
    $k([\psi];A)\sg_N=0$, $k([\psi],A)\in {\bR}$;

\item[ii)] if $H^1(Q^n;{\bR})=0$ and $[\psi]\neq 0$, then $\sg_N=0$.
\end{enumerate}
\end{th}

%   PROOF
\begin{pf}
In condition i) we have $[\psi]=k([\psi];A)\sg_A$ so that
$k([\psi];A)\sg_N=-G^* \wh{i}^*[\psi]=-[\bt]=0$.

If $[\psi]\neq 0$, then $k([\psi];A)\neq 0$ and ii) follows from i).
\end{pf}

%   REMARK 2.9
\begin{rem}
When $M={\bR}^{2n}$, all conditions of Theorem 2.5 are satisfied.
The form $\psi$ is a closed 1-form, more precisely, $\psi=d\tht$
(see formula 2.8), where $\tht$ is the function from
$L(T{\bR}^{2n})={\bR}^{2n}\times U(n)/O(n)$ to $S^1$. It is not
difficult to calculate that in this case $k([\psi];A)=2\pi$.
Therefore Theorem 3 in [8] (Fomenko's conjecture) follows from
Theorem 2.8 as a corollary not only for minimal Lagrangian
submanifolds in ${\bR}^{2n}$ as in [8], but already for the
manifolds from $LH(M)$ as well.
\end{rem}

By Remark 2.9 and Theorem 2.8 we are able to extend Fomanko's conjecture
for the manifolds from the class $LH(M^{2n})$. When the form $\psi$
is a closed 1-form, this conjecture equivalently can be reduced to the
exactness of the form $A^*[\psi]|_N$.

%   COROLLARY 2.10
\begin{cor}
If $M^n$ is a Kaehlerian manifold, with $H^1(M^n;{\bR})=0$, and
Ricchi's tensor is identically equal to zero, then for $N\in
LH(TM^n)$ we have $k([\psi],A)\sg_N=0$. In the particular case,
where $[\psi]\neq 0$, we have $\sg_N=0$.
\end{cor}

This corollary follows from Theorems 2.4 and 2.8.

\section*{Aknowledgement}

The research described in this  paper was made possible in part by
Grant No. BJV000 from the International Science Foundation.
\vskip+0.8cm

\centerline{\sc References}
\vskip+0.4cm

1. V. P. Maslov, Theory of perturbation and asymptotic methods. (Russian)
    {\em Moscow University Press, Moscow,} 1965.

2. V. I. Arnold, Characteristic class entering in quantization conditions.
    (Russian) {\em Funkts. Anal. Prilozh.} {\bf 1}(1967), No. 1, 1--14.

3. D. B. Fuks, Maslov--Arnold characteristic classes. (Russian)
    {\em Dokl. Akad. Nauk SSSR} {\bf 178}(1968), No. 2, 303--306.

4. V. V. Trofimov, Symplectic connections, index of Maslov and conjecture
    of Fomenko. (Russian)
    {\em Dokl. Akad. Nauk SSSR} {\bf 304}(1989), No. 6, 1302--1305.

5. L. H\"{o}rmander, Fourier integral operators I.
    {\em Acta Math. } {\bf 127}(1971), 79--183.

6. L. H\"{o}rmander, The calculus of Fourier integral operators.
    {\em Prospects in Math.,} 33--57, {\em Ann. of Math. Stud.,} 70,
    {\em Princeton Univ. Press., Princeton, N.~Y., } 1971.

7. V. Guilleman and S. Sternberg, Geometric Asymptotics.
    {\em Amer. Math. Soc., Providence, Rhode Island,} 1977.

8. Le Hong Van and A. T. Fomenko, The criteria of minimality
    of Lagrangian submanifolds in a Kaehlerian manifold. (Russian)
    {\em Mat. Zametki} {\bf 4}(1987), 559--571.

9. Z. A. Tevdoradze, One note on the calculation of Arnold--Maslov
    class of Liouville tori in cotangent bundles. (Russian)
    {\em Uspekhi Mat. Nauk} {\bf 43}(1988), No. 5, 219--220.

10. Z. A. Tevdoradze, Calculation of Maslov class in the Lagrange and
    Kovalevskaya problem for the motion of a rigid body about
    a fixed point. (Russian) {\em Bull. Acad. Sci. Georgia}
    {\bf 141}(1991), No. 1, 45--49.

11. E. Dyer, Cohomology theories. {\em W. A. Benjamin Inc., New York,} 1969.

12. A. Lichnerowicz, Th\'{e}orie globale des connexions et des groups
    d'ho\-lo\-no\-mie. {\em Ed. Cremonese, Roma,} 1955.

13. Le Hong Van, Minimal $\Phi$-Lagrangian surfaces in almost Hermitian
    manifolds. (Russian) {\em Mat. Sb.} {\bf 180}(1989), 924--936.

14. R. Harvey and H. B. Lawson, Calibrated geometries.
    {\em Acta Math.} {\bf 148}(1982), 47--157.

\

%\centerline{(Received 28.11.1994; revised version 15.03.1995)}

\

\end{document}